\documentclass[11pt]
{article}
\usepackage{amsthm,amsfonts,amssymb,latexsym}

\newcommand{\blankline}{\mbox{}\hfill\\ \mbox{}}
\newcommand{\n}{\noindent}
\newcommand{\vn}{\vspace{5mm} \noindent}
\def\RR{{\mathbb R}}

\def\QQ{{\mathbb Q}}

\def\ZZ{{\mathbb Z}}

\def\GG{{\mathbb G}}

\def\DD{{\mathbb D}}
\def\FF{{\mathbb F}}

\def\BT1{{$\mbox{\rm BT}_{1}$}}
\def\DM1{{$\mbox{\rm DM}_{1}$}}

\newcommand{\Spec}{\mbox{\rm Spec}}
\newcommand{\End}{\mbox{\rm End}}

\newcommand{\Ker}{\mbox{\rm Ker}}
\newcommand{\Image}{\mbox{\rm Im}}

\newcommand{\mini}{\mbox{\rm min}}

\newcommand{\va}{\varphi}

\def\cF{{\cal F}}
\def\cV{{\cal V}}

\def\cN{{\cal N}}

\newcommand{\B}{{\hfill$\Box$}}   

\makeatletter
\renewcommand{\subsection}[1]{\@startsection{subsection}{2}{\z@}
            {-3.25ex plus -1ex minus -.2ex}{-1sp}{\normalsize\bf}
             {\ignorespaces#1 }}
 
\makeatother
\begin{document}
\sloppy

\title{Minimal p-divisible groups}
\author{Frans Oort}
\date{Version 16-III-2004}
\maketitle

\section*{}{\bf Introduction.}   A   $p$-divisible group   $X$ can be seen as a tower of building blocks, each of which is isomorphic to the same finite group scheme $X[p]$. Clearly, if $X_1$ and $X_2$ are isomorphic then $X_1[p] \cong X_2[p]$; however, conversely $X_1[p] \cong X_2[p]$ does  in general  not imply that $X_1$ and $X_2$ are isomorphic. Can we give, over an algebraically closed field in characteristic    $p$, a condition  on the $p$-kernels which ensures this converse? Here are two known examples of such a condition: consider the case that $X$ is {\it ordinary}, or the case that $X$ is {\it superspecial} ($X$ is the $p$-divisible group of a product of supersingular elliptic curves); in these cases the
$p$-kernel  uniquely determines $X$. 

These are special cases of  a surprisingly  complete and simple answer: 
\begin{center}
{\it if $G$ is ``minimal'', then $X_1[p] \cong G \cong X_2[p]$ implies   $X_1 \cong X_2$,}
\end{center}
see \ref{th}; for a definition of ``minimal'' see \ref{notat}. This is ``necessary and sufficient'' in the  sense that any $G$ that is  {\it not minimal} there exist infinitely many mutually non-isomorphic   $p$-divisible groups   with    $p$-kernel isomorphic to $G$; see \ref{Rk41}.

\vn 
{\bf Remark} (motivation). You might wonder why this is interesting. 
\begin{description}
\item[{\bf EO}]  In \cite{EO} we have defined a natural {\it stratification} of the moduli space of polarized abelian varieties in positive characteristic: moduli points are in the same stratum if and only if the corresponding $p$-kernels are geometrically isomorphic. Such strata are called EO-strata. 
\item[{\bf Fol}] In \cite{FO-Foliations} we  define in the same moduli spaces a {\it foliation} : moduli points are in the same leaf if and only if the corresponding $p$-divisible groups are geometrically isomorphic; in this way we obtain a foliation of every open Newton  polygon stratum. 
\item[Fol $\subset$ EO] The observation $X \cong Y \Rightarrow
X[p] \cong Y[p]$ shows that
any leaf in the second sense is contained in precisely one stratum in the first sense; the main result of this paper,  ``$X$ {\it is minimal if and only if $X[p]$ is minimal}'', shows that {\it a stratum} (in the first sense) {\it and a leaf} (in the second  sense) {\it are equal} if we are in the minimal, principally polarized situation.
\end{description} 

\vn
In this paper  we consider    $p$-divisible groups   and finite group schemes over an {\it algebraically closed } field $k$ of characteristic    $p$.

\vn
An apology. In \ref{notat3} and in \ref{notat4} we fix notations, used for the proof of \ref{filt}, respectively \ref{ext}; according to the need, the notations in these two different cases are  different. We hope this difference in notations in Section \ref{SF} versus Section \ref{EP} will not cause confusion. 

\vn
 Group schemes considered are supposed to be commutative. We use {\it covariant} Dieudonn\'e module theory. We write $W = W_{\infty}(k)$ for the ring of infinite Witt vectors with coordinates in $k$. Finite products in the category of $W$-modules are denoted ``$\times$'' or by ``$\prod$'', while  finite products in the category of Dieudonn\'e modules are denoted by ``$\oplus$''; for finite products of $p$-divisible groups we use ``$\times$'' or  ``$\prod$''. We write $F$ and $V$, as usual,  for ``Frobenius'' and ``Verschiebung'' on commutative group schemes; we write  $\cF = \DD(V)$ and $\cV = \DD(F)$, see \cite{EO}, 15.3,  for the corresponding operations on Dieudonn\'e modules.

\vn
{\bf Acknowledgments.} {\footnotesize Part of the work for this paper was done while visiting  Universit\'e de Rennes, and the Massachusetts Institute of Technology; I thank the Mathematics Departments of these universities for hospitality and stimulating working environment. I thank Bas Edixhoven and Johan de Jong for discussions on ideas necessary for this paper.}

\section{Notations and the main result.}

\subsection{}\label{notat}{\bf Some definitions and notations.}\\
$H_{m,n}$. \ \ \  
We define the   $p$-divisible group   $H_{m,n}$ over the prime field $\FF_p$ in case $m$ and $n$ are coprime non-negative integers. In case $m$ and $n$ are positive this is given by $\DD(H_{m,n}) := E/E(\cF^m - \cV^n)$, see \cite{AJdJ.FO}, 5.2;  here $E$ is  the      Dieudonn\'e ring, e.g. see \cite{Manin}, Th. 1.4 on pp. 21/22. This $p$-divisible group   $H_{m,n}$ is of dimension $m$, its Serre-dual $X^t$ is of dimension $n$, it is isosimple, and {\it its endomorphism ring $\End(H_{m,n} \otimes \overline{\FF_p})$ is the maximal order in the endomorphism algebra} $\End^0(H_{m,n} \otimes \overline{\FF_p}$) (and these properties characterize this   $p$-divisible group   over  $\overline{\FF_p}$). We will use the notation $H_{m,n}$ over any base $S$ in characteristic    $p$, i.e. we write $H_{m,n}$ instead of $H_{m,n} \times_{\Spec(\FF_p)} S$, if no confusion can occur.

The ring $\End(H_{m,n} \otimes \FF_p) = R'$ is commutative; write $L$ for the field of fractions of $R'$. Consider integers $x, y$ such that for the coprime positive integers $m$ and $n$ we have $x{\cdot}m + y{\cdot}n = 1$. In $L$ 
we define the element $\pi = \cF^y{\cdot}\cV^x \in L$. Write $h = m+n$. Note that
$\pi^h = p$ in $L$. Here $R' \subset L$ is the maximal order, hence $R'$ integrally closed in $L$, and we conclude that $\pi \in R'$.  This element $\pi$ will be called the uniformizer in  this endomorphism ring. In fact, $W_{\infty}(\FF_p) = \ZZ_p$, and $R' \cong \ZZ_p[\pi]$ and $L \cong  \QQ_p[[\cF,\cV]]/(\cF^m-\cV^n, \cF\cV-p)$.  In $L$ we have:
$$m+n =: h, \ \ \pi^h = p, \ \ \cF = \pi^n, \ \ \cV = \pi^m.$$
For a further description of $\pi$, of $R = \End(H_{m,n} \otimes k)$ and of  $D = \End^0(H_{m,n} \otimes k)$ see \cite{AJdJ.FO}, 5.4; note that $\End^0(H_{m,n} \otimes k)$ is non-commutative if $m>0$ and $n>0$. Note that $R$ is a ``discrete valuation ring'' (terminology sometimes also used for non-commutative rings).

\vn
{\bf Newton polygons.}\ \ Let $\beta$ be a   Newton polygon. By definition, in the notation used here, this is a lower convex polygon in $\RR^2$ starting at $(0,0)$,  ending at $(h,c)$ and having break points with integral coordinates; it is given by $h$ slopes in non-decreasing order; every slope $\lambda$ is a rational number, $0 \leq \lambda \leq 1$. 

To each ordered pair of nonnegative integers $(m,n)$ we assign a set of  $m+n=h$ slopes equal to $n/(m+n)$; this Newton polygon ends at $(h,c=n)$. 

In this way a   Newton polygon corresponds with  a set of pairs $\sum_i \   (m_i,n_i)$ and conversely. Usually we consider only coprime  pairs $(m_i,n_i)$;  we write $H(\beta) := \times_i \   H_{m_i,n_i}$ in  case $\beta = \sum_i \   (m_i,n_i)$. A $p$-divisible group $X$ over a field of positive characteristic defines a  Newton polygon where $h$ is the height of $X$ and $c$ is the dimension of its Serre-dual $X^t$. By the Dieudonn\'e-Manin classification:  two $p$-divisible groups over an algebraically closed field of positive characteristic are isogenous if and only if their Newton polygons are equal.

\vn
{\bf Definition.} {\it A   $p$-divisible group   $X$ is called}  minimal {\it  if there exists a   Newton polygon $\beta$ and an isomorphism $X_k \cong H(\beta)_k$, where $k$ is an algebraically field.} 

Note that in every isogeny class of $p$-divisible groups over an algebraically closed field there is precisely one minimal $p$-divisible group.

\vn
{\bf Truncated $p$-divisible groups.}  A finite group scheme $G$ (finite and flat over some base, but in this paper we will soon work over a field) is called a      $\mbox{\rm BT}_{1}$, see \cite{Ill}, page 152, if $G[F] := \Ker F_G = \Image V_G =: V(G)$  and   $G[V] = F(G)$ (in particular this implies that $G$ is annihilated by    $p$). Such group schemes over a perfect field appear as the    $p$-kernel of a   $p$-divisible group, see \cite{Ill}, Prop. 1.7 on page 155. The abbreviation ``$\mbox{\rm BT}_{1}$''  stand for ``1-truncated Barsotti-Tate group'';  the terms ``$p$-divisible group'' and ``Barsotti-Tate group'' indicate the same concept.

The      Dieudonn\'e module of a     $\mbox{\rm BT}_{1}$ over a perfect field $K$ is called a     $\mbox{\rm DM}_{1}$; for $G = X[p]$ we have $\DD(G) = \DD(X)/p\DD(X)$. In other terms: such a      Dieudonn\'e module $M_1 = \DD(X[p])$ is a finite dimensional vector space over $K$, on which $\cF$ and $\cV$ operate in the usual way, with the property that $M_1[\cV] = \cF(M_1)$  and $M_1[\cF] = \cV(M_1)$.

\vn
{\bf Definition.}  {\it A     $\mbox{\rm BT}_{1}$ \  $G$ is called}  minimal {\it if there exists a   Newton polygon  $\beta$ such that $G_k \cong H(\beta)[p]_k$. A     $\mbox{\rm DM}_{1}$ is called} minimal {\it if it is the     Dieudonn\'e module of a minimal     $\mbox{\rm BT}_{1}$.}

\subsection{}{\bf Theorem.}\label{th} {\it Let $X$  be    a $p$-divisible group over an algebraically closed field $k$ of characteristic    $p$. Let $\beta$ be 
a Newton polygon. Then
$$X[p] \cong H(\beta)[p] \quad\Longrightarrow\quad X \cong  H(\beta).$$
In particular: if $X_1$ and $X_2$ are  $p$-divisible groups over $k$, with $X_1[p] \cong G \cong X_2[p]$, where $G$ is }  minimal, {\it then $X_1 \cong X_2$.}\\
{\bf Remark.} We have no a priori condition on the Newton polygon of $X$, nor do we a priori assume that $X_1$ and $X_2$ have the same Newton polygon.\\
{\bf Remark.} In general an isomorphism $\va_1: X[p] \to H(\beta)[p]$ does not lift to an isomorphism $\va: X \to H(\beta)$.

\subsection{} Here is another way of explaining the result of this paper. Consider the map
$$[p]: \{X \mid \mbox{a $p$-divisible group}\}/\cong_k  \quad\longrightarrow \{G \mid \mbox{a \BT1}\}/\cong_k,   \quad \quad X \mapsto X[p].$$
This map is surjective, e.g. see \cite{Ill}, 1.7; also see \cite{EO}, 9.10. 
\begin{itemize}
\item By results of this paper we know:  {\it For every Newton polygon $\beta$ there is an isomorphism class $X := H(\beta)$ such that the fiber of the map  $[p]$ containing $X$ consists of one element.} 
\item {\it For every $X$ not isomorphic to some $H(\beta)$ the fiber of $[p]$ containing $X$ is infinite;} see \ref{Rk41} 
\end{itemize}

\vn
{\bf Convention.} The slope $\lambda= 0$, given by the pair $(1,0)$, defines the $p$-divisible group $G_{1,0} = \GG_m[p^{\infty}]$, and its $p$-kernel is $\mu_p$. The slope $\lambda = 1$, given by the pair $(0,1)$, defines the $p$-divisible group $G_{0,1}  = \underline{\QQ_p/\ZZ_p}$ and its $p$-kernel is $\underline{\ZZ/p}$. These $p$-divisible groups and their $p$-kernels split off naturally over a perfect field, see \cite{CGS}, 2.14. The theorem is obvious for these minimal \BT1 \  group schemes over an algebraically closed field. It suffices to prove the theorem in case all group schemes considered are of local-local type, i.e. all slopes considered are strictly between $0$ and $1$; from now on we make these assumptions.

\subsection{}\label{A} We give already one explanation about notation and method of proof.  Let $m, n \in \ZZ_{>0}$ be coprime. Start with 
$H_{m,n}$ over $\FF_p$. Let $Q' = \DD(H_{m,n} \otimes \FF_p)$. In the terminology of \cite{AJdJ.FO}, 5.6 and Section 6, a semi-module of $H_{m,n}$ equals $[0,\infty) = \ZZ_{\geq 0}$. Choose a non-zero element in $Q'/\pi Q'$, this is a on-dimensional vector space over $\FF_p$, and lift this element to $A_0 \in Q'$. Write $A_i = \pi^i A_0$ for every $i \in \ZZ_{>0}$. Note that
$$\pi A_i = A_{i+1}, \ \ \  \cF A_i = A_{i+n}, \ \ \  \cV A_i = A_{i+m}.$$  
Fix an algebraically closed field $k$; we write $Q = \DD(H_{m,n} \otimes k)$. Clearly $A_i \in Q' \subset Q$, and the same relations as given above hold.
Note that $\{A_i \mid i \in \ZZ_{\geq 0}\}$   generate $Q$ as a $W$-module. {\it The fact that a semi-module of the minimal $p$-divisible group $H_{m,n}$ does not contain ``gaps''  is the essential (but sometimes hidden) argument in the proofs below.} 

The set $\{A_0, \cdots , A_{m+n-1}\}$ is a $W$-basis for $Q$. If $m \geq n$ we see that $\{A_0, \cdots , A_{n-1}\}$ is a set of generators for $Q$ as a  Diedonn\'e module;  the structure of this Diedonn\'e module can be described as follows; for this set of generators   we consider  another numbering $\{C_1, \cdots , C_n\} = \{A_0, \cdots , A_{n-1}\}$ and we define positive integers $\gamma_i$ by: $C_1 = A_0$ and $\cF^{\gamma_1} C_1 = \cV C_2, \cdots , \cF^{\gamma_n} C_n = \cV C_1$ (note that we assume $m \geq n$),  which gives a ``cyclic'' set of generators for $Q/pQ$ in the sense of \cite{Kraft}. These notations will be repeated and explained more in detail in \ref{notat3} and \ref{notat4}.

\section{A slope filtration}\label{SF}
\subsection{}\label{notat2} We consider a Newton polygon $\beta$ given by $r_1(m_1,n_1), \cdots , r_t(m_t,n_t)$; here $r_1, \cdots , r_t \in \ZZ_{>0}$, and every $(m_j,n_j)$ is an ordered pair of coprime positive integers; we write $h_j = m_j +n_j$ and we suppose the ordering is chosen in such a way that $\lambda_1 := n_1/h_1 < \cdots < \lambda_t := n_t/h_t$. Write 
$$H := H(\beta) = \prod_{1 \leq j \leq t} \  (H_{m_j,n_j})^{r_j};  \ \ G := H(\beta)[p].$$ 
The following proposition uses this notation; suppose that $t>0$.

\subsection{}\label{filt}{\bf Proposition.}  {\it Suppose $X$ is a $p$-divisible group over an algebraically closed field $k$. Suppose that $X[p] \cong H(\beta)[p]$. Suppose that $\lambda_1 = n_1/h_1 \leq 1/2$. Then there exists a $p$-divisible subgroup $X_1 \subset X$ and  isomorphisms 
$$X_1 \cong (H_{m_1,n_1})^{r_1} \quad\mbox{and}\quad (X/X_1)[p] \cong \prod_{j>1}\ (H_{m_j,n_j}[p])^{r_j}.$$}

\subsection{}{\bf Remark.} The condition that $X[p]$ is minimal is essential; e.g. it is easy to give an example of a $p$-divisible group $X$ which is isosimple, such that $X[p]$ is decomposable.

\subsection{}\label{filt2}{\bf Corollary.} {\it For $X$  with $X[p] \cong H(\beta)[p]$, with $\beta$  as in} \ref{notat2}, {\it there exists a filtration by $p$-divisible subgroups
$$X_0 := 0 \subset X_1 \subset \cdots \subset X_t = X \quad\mbox{such that}\quad X_{j}/X_{j-1} \cong (H_{m_j,n_j})^{r_j}, \quad\mbox{for}\quad  1 \leq j \leq t.$$}
{\bf Proof of the corollary.} Assume by induction that the result has been proved for all $p$-divisible groups where $Y[p] = H(\beta')[p]$ is minimal such that $\beta'$ has at most $t-1$ different slopes; induction starting at $t-1=0$, i.e. $Y = 0$. If on the one hand the smallest slope of $X$ is at most $1/2$, the proposition gives $0 \subset X_1 \subset X$, and using the induction hypothesis on $Y = X/X_1$ we derive the desired filtration. If  on the other hand all slopes of $X$ are bigger than $1/2$, we apply the proposition to the Serre-dual of $X$, using the fact that the Serre-dual of $H_{m,n}$ is $H_{n,m}$;  dualizing back we obtain $0 \subset X_{t-1} \subset X$, and using the induction hypothesis on $Y = X_{t-1}$ we derive the desired filtration. Hence we see that the proposition gives the induction step; this proves the corollary.  \B\ref{filt}$\Rightarrow$\ref{filt2} 

\subsection{}\label{notat3} We use notation as in \ref{notat2} and \ref{filt}, and we fix further notation which will be used in the proof of \ref{filt}.   

\vn
{\it Let $M = \DD(X)$.  We write $Q_j =\DD(H_{m_j,n_j})$. Hence 
$$M/pM \cong \bigoplus_{1\leq j\leq t} \ \ (Q_j/pQ_j)^{r_j}.$$
Using this isomorphism we construct a map} 
$$v: M  \longrightarrow \QQ_{\geq 0} \cup\{\infty\}.$$

\vn
Let $\pi_j$ be the uniformizer of $\End(Q_j)$, see \ref{notat}. As in \ref{A} we choose $A^{(j)}_i \in  Q_j$ with $i \in \ZZ_{\geq 0}$
which generate $Q_j$ such that $\pi_j{\cdot}A^{(j)}_i = A^{(j)}_{i+1}$, \ \ $\cF{\cdot}A^{(j)}_i =A^{(j)}_{i+n_j}$  and $\cV{\cdot}A^{(j)}_i=A^{(j)}_{i+m_j}$ . We have $Q_j/pQ_j = \times_{0 \leq i < h_j}  \ k{\cdot}(A^{(j)}_i \ {\rm mod} \ pQ_j)$. 

\vn
We write 
$$A^{(j)}_i = (A^{(j)}_{i,s} \mid 1 \leq s \leq r_j) \in (Q_j)^{r_j}$$
for the vector with coordinate $A^{(j)}_{i,s}$ in the factor $s$.
For  $B \in M$   we uniquely write 
$$B \ {\rm mod} \ pM = a = \sum_{j, \ 0 \leq i< h_j, \ 1 \leq s \leq r_j} \  b^{(j)}_{i,s}{\cdot}(A^{(j)}_{i,s} \ {\rm mod} \  pQ_j), \quad b^{(j)}_{i,s} \in k;$$
if moreover $B \not\in pM$ we define
$$v(B) = {\rm min}_{j,\ i,\ s, \  b^{(j)}_{i,s} \not= 0}\ \ \frac{i}{h_j}.$$ If $B' \in p^{\beta}M$ and $B' \not\in p^{\beta + 1}M$ we define $v(B') = \beta + v(p^{-\beta}{\cdot}B')$. We write $v(0) = \infty$. {\it This ends the construction of} $v: M  \longrightarrow \QQ_{\geq  0} \cup\{\infty\}$.

\vn  
For any $\rho \in \QQ$ we define 
$$M_{\rho} = \{B \mid v(B) \geq \rho\};$$
note that $pM_{\rho} \subset M_{\rho + 1}$.
Let $T$ be the least common multiple of $h_1, \cdots ,\ h_t$. Note that, in fact, $v: M - \{0\} \to \frac{1}{T} \ZZ_{\geq  0}$. Note that, by construction, $v(B) \geq d \in \ZZ$ if and only if $p^d$ divides $B$ in $M$. Hence $\cap_{\rho \to \infty} \ \ M_{\rho}= \{0\}$. 

\vn
The basic assumption $X[p] \cong H(\beta)[p]$ of \ref{th} is:
$$M/pM = \bigoplus_{1 \leq j \leq t,  \ 1 \leq  s \leq r_j} \ \ \  \prod_{0 \leq i <h_j} \ \  k{\cdot}((A^{(j)}_{i,s} \ {\rm mod} \  pQ_j))$$
(we write this isomorphism of Dieudonn\'e modules as an equality).
 For $ 0 \leq i < h_j$ and $1 \leq s \leq r_j$  we choose $B^{(j)}_{i,s} \in M$ such that:  
$$B^{(j)}_{i,s}  \ {\rm mod} \ pM = A^{(j)}_{i,s} \ {\rm mod} \  pQ_j.$$ 
Define  $B^{(j)}_{i +\beta{\cdot}h_j,s} = p^{\beta}{\cdot}B^{(j)}_{i,s}$. By construction we have:  $v(B^{(j)}_{i,s}) = i/h_{j}$ for all $i \geq 0$, all $j$ and all $s$. Note that $M_{\rho}$ is generated over $W = W_{\infty}(k)$ by all elements $B^{(j)}_{i,s}$ with $v(B^{(j)}_{i,s}) \geq \rho$. As a short-hand we will write $$B^{(j)}_i \textrm{  for the vector  } (B^{(j)}_{i,s}\mid 1 \leq s \leq r_j) \in M^{r_j}.$$

\vn
We write $P \subset M$ for the sub-$W$-module generated by all $B^{(j)}_{i,s}$ with $j \geq 2$ and $i < h_j$; we write $N \subset M$  for the sub-$W$-module generated by all $B^{(1)}_{i,s}$ with $i< h_1$. Note that
$M = N \times P$, a direct sum of $W$-modules. Note that $M_{\rho} = (N \cap M_{\rho}) \times (P \cap M_{\rho})$.

\vn
{\it In the proof the $W$-submodule $P \subset M$ will be fixed;   its $W$-complement $N \subset M$ will change eventually if it is not already a Dieudonn\'e submodule.}

\vn
We write $m_1 = m, \ \ n_1 = n,\ \  h =  h_1 = m + n$, and $r = r_1$. Note that we assumed $0 < \lambda_1 \leq 1/2$, hence $m \geq n > 0$.
For $i \geq 0$ we  define integers $\delta_i$ by:
$$i{\cdot}h \leq \delta_i{\cdot}n < i{\cdot}m + (i+1){\cdot}n = ih + n$$
and non-negative integers $\gamma_i$ such that 
$$\delta_1 = \gamma_1 + 1, \cdots ,\delta_i = \gamma_1 + 1 + \gamma_2 + 1 + \cdots +\gamma_i + 1, \cdots ;$$
note that $\delta_n = h = m+n$; hence $\gamma_1 + \cdots + \gamma_n = m$.  For $1 \leq i \leq n$ we write 
$$f(i) = \delta_{i-1}{\cdot}n - (i-1){\cdot}h;$$ this means that $0 \leq f(i) < n$ is the remainder of dividing $\delta_{i-1}n$ by $h$; not that $f(1) = 0$. As $\gcd(n,h) =1$ we see that  
$$f :\{1, \cdots , n\} \to \{0, \cdots ,n-1\}$$ 
is a bijective map. The inverse map $f'$ is given by:  $$f': \{0, \cdots ,n-1\} \to \{1, \cdots , n\}, \quad  f'(x) \equiv 1 - \frac{x}{h} \pmod{n}, \quad 1 \leq f'(x) \leq n.$$

\vn
In $(Q_1)^r$ we have the vectors $A^{(1)}_i$.  
We choose  $C'_1:=A^{(1)}_0$ and we choose  $\{C'_1, \cdots ,C'_n\} = \{A^{(1)}_0, \cdots ,A^{(1)}_{n-1}\}$ by 
$$C'_i := A^{(1)}_{f(i)}, \quad C'_{f'(x)} = A^{(1)}_x;$$ this means that:
$$\cF^{\gamma_i} C'_i = \cV C'_{i+1}, \ \ 1 \leq i < n, \ \ \ \cF^{\gamma_n}C'_n = \cV C'_1,  \quad\mbox{hence}\quad \cF^{\delta_i} C'_1 = p^i{\cdot}C'_{i+1}, \ \ \ 1 \leq i < n;$$
 note that
$\cF^h C'_1 = p^n{\cdot} C'_1$. With these choices we see that 
$$\{\cF^j C'_i \mid 1 \leq i \leq n, 0 \leq j \leq \gamma_i\} = \{A^{(1)}_{\ell} \mid 0 \leq \ell < h\}.$$
For later reference we state:

\subsection{}\label{struct} {\it Suppose $Q$ is a Dieudonn\'e module  with an element $C \in Q$, such that there exist coprime integers $n$ and $n+m = h$  as above such that $\cF^h{\cdot}C = p^n{\cdot}C$ and such that $Q$ as a $W$-module is generated by $\{p^{-[jn/h]}\cF^jC  \mid 0 \leq j < h\}$, then $Q \cong \DD(H_{m,n})$.}\\
This is proved by explicitly writing out the required isomorphism. Note that $\cF_n$ is injective on $Q$, hence $\cF^h{\cdot}C = p^n{\cdot}C$ implies $\cF^m{\cdot}C = \cV^n{\cdot}C$.

\subsection{} Accordingly we choose  $B^{(1)}_{f(i),s} =: C_{i,s} \in M$ with $1 \leq i \leq n$. 
Note that
$$\{\cF^j C_{i,s} \mid 1 \leq i \leq n, \ \ 0 \leq j \leq \gamma_i \ \ 1 \leq s \leq r\} \quad\mbox{is a}\quad W \mbox{-basis for}\quad N,$$
$$\cF^{\gamma_i}C_{i,s} - \cV C_{i+1,s} \in pM, \ \ 1 \leq i<n, \ \
\cF^{\gamma_n}C_{n,s} - \cV C_{1,s} \in pM.$$
We write $C_i = (C_{i,s} \mid 1 \leq s \leq r)$. 
As a reminder, we sum up some of the notation constructed:
$$
\begin{array}{rcc}
N \subset M \quad\ \  &  & \bigoplus_j \ \ (Q_j)^{r_j}   \\
\Big\downarrow \quad\ \   &   &  \Big\downarrow  \\
M/pM  &  =  &    \bigoplus_j \ \ (Q_j/pQ_j)^{r_j}, \\
B^{(j)}_{i,s}  \in M &    & A^{(j)}_{i,s} \in Q_j\\
 C_{i,s} \in N & &  C'_{i,s} \in Q_1.
\end{array}
$$


\subsection{}{\bf Lemma.}\label{L1} {\it Use the notation fixed up to now.}\\
{\bf (1)} {\it  For every $\rho \in \QQ_{\geq  0}$  the map $p: M_{\rho} \to M_{\rho + 1}$, multiplication by $p$, is surjective. }\\
{\bf (2)} {\it  For every $\rho \in \QQ_{\geq  0}$ we have $\cF M_{\rho} \subset M_{\rho + (n/h)}$.}\\
{\bf (3)} {\it For every $i$ and $s$ we have $\cF B^{(1)}_{i,s} \in M_{(i+n)/h}$;  for every $i$ and $s$ and every $j > 1$  we have 
$\cF B^{(j)}_{i,s} \in M_{(i/h_j) + (n/h) + (1/T)}$.}\\
{\bf (4)} {\it   For every $1 \leq i \leq  n$   we have
$\cF^{\delta_i} C_1 - p^i B^{(1)}_{f(i+1)} \in
(M_{i + (1/T)})^r$; moreover  $\cF^{\delta_n} C_1 - p^n C_1 \in
(M_{n + (1/T)})^r$.}\\
{\bf (5)} {\it If $u$ is an integer with $u>Tn$, and $\xi_N \in (N \cap M_{u/T})^r$, there exists 
$$\eta_N \in N \cap (M_{(u/T)-n})^r \quad\mbox{such that}\quad 
(F^h - p^n)\eta_N \equiv \xi_N \pmod{(M_{(u+1)/T})^r}.$$  
} \\
{\bf Proof.} We know that   $M_{\rho + 1}$ is generated by the elements $B^{(j)}_{i,s}$ with $i/h_j \geq \rho + 1$; because $\rho \geq 0$ such elements satisfy $i \geq h_j$. Note that
$p{\cdot}B^{(j)}_{i-h_j,s} =  B^{(j)}_{i,s}$. This proves the first property.

\B(1)

\vn
At first we show $\cF M \subset M_{n/h}$. Note that 
for all $1 \leq j \leq t$ and all $\beta \in \ZZ_{\geq 0}$
$$\beta h_j \leq i < \beta h_j + m_j \quad\Rightarrow\quad
\cF B^{(j)}_i = B^{(j)}_{i+n_j},  \hspace{55mm} (\ast)
$$
and 
$$\beta h_j + m_j \leq i < (\beta + 1) h_j \quad\Rightarrow\quad B^{(j)}_i = \cV B^{(j)}_{i-m_j} + p^{(\beta + 1)}\xi, \ \ \ \xi \in M^{r_j}.  \hspace{15mm}(\ast\ast)$$
from these properties, using $n/h \leq n_j/h_j$ we conclude: \ $\cF M \subset M_{n/h}$. 

Further we see:
by $(\ast)$ we have
$$v(\cF B^{(j)}_{i,s}) = v(B^{(j)}_{i+n_j,s}) = (i+n_j)/h_j,$$
and  
$$\frac{i+n_j}{h_j}= \frac{i+n}{h}  \quad\mbox{if}\quad j=1; \quad \      \frac{i+n_j}{h_j} > \frac{i}{h_j} + \frac{n}{h}  
\quad\mbox{if}\quad j>1.$$
By $(\ast\ast)$ it suffices to consider only $m_j \leq i < h_j$, and  hence $\cF B^{(j)}_{i,s} = p  B^{(j)}_{i-m_j,s} + p \cF \xi$; so we have
$$v(\cF B^{(j)}_{i,s}) \geq \mini\left(v(p  B^{(j)}_{i-m_j,s}), v(p\cF \xi_s) \right);$$
for $j=1$ we have $v(p  B^{(1)}_{i-m_1,s}) = (i+n)/h \geq 1$ and $v(p \cF \xi) \geq 1 + (n/h) > (i/h) + (n/h)$; for $j>1$ we have $v(p  B^{(j)}_{i-m_j,s}) > (i/h_j) + (n/h)$ and $(i/h_j) + (n/h) < 1 + (n/h) \leq v(p\cF \xi_s)$; hence $v(\cF B^{(j)}_{i,s})) > (i/h_j) + (n/h)$ if $j > 1$. This ends the proof of (3). Using (3) we see that (2) follows.  \B(2)+(3)

\vn
{}From $\cF^{\gamma_i} C_i =  \cV C_{i+1} + \xi_i$ for $i<n$ and  $\cF^{\gamma_i} C_n =  \cV C_1 + \xi_n$, here $\xi_i \in M^r$ for $i \leq n$,   we conclude:
$$\cF^{\delta_i} C_1 = p^i C_{i+1} + \sum{}_{1 \leq \ell \leq i} \ \ p^{\ell} \cF^{\delta_i - \delta _{\ell}} \cF \xi_{\ell},$$
and the analogous fomula for $i = n$ (write $C_{n+1}=C_1$). Note that 
$$ih \leq \delta_i n \quad\mbox{and}\quad \delta_{\ell} n < \ell m + (\ell+1)n = \ell h + n;$$ 
this shows that
$$\ell h + (\delta_i - \delta_{\ell})n + n > ih;$$
using (2) we conclude (4). \B(4)

\vn
Note that $h = h_1$ divides $T$. If $\ell$ is an integer such that $(\ell - 1)/h < u/T < l/h$ then $u < u+1 \leq \ell \frac{T}{h}$; in this case we see that $N \cap M_{u/T} = N \cap M_{(u+1)/T}$. In this case we choose $\eta_N = 0$.

Suppose that $\ell$ is an integer with $u/T = \ell/h$. Then $N \cap M_{u/T} = N_{\ell/h} \supset N_{(\ell+1)/h} = N \cap M_{(u+1)/T}$. We consider the image of $N \cap M_{(\ell/h)-n}$ under $F^h - p^n$. We see, using previous results, that this image is in $N_{\ell/h} + M_{(u+1)/T}$ (here ``+'' stands for the span as $W$-modules). We obtain a factorization and an isomorphism
$$F^h - p^n: N \cap M_{(\ell/h)-n} \longrightarrow 
\left(N_{\ell/h} + M_{(u+1)/T}\right) / M_{(u+1)/T} \cong N_{\ell/h}/N_{(\ell+1)/h}.$$
{\it We claim that this map is surjective.}  The factor space $N_{\ell/h}/N_{(\ell+1)/h}$ is a vector space over $k$ spanned by the residue classes of the  elements $B^{(1)}_{\ell,s}$. For the residue class of $y_s B^{(1)}_{\ell,s}$ we solve the equation $x_s^{p^n}- x_s = y_s$ in $k$; lifting these $x_s$ to $W$ (denoting the lifts by the same symbol), we see that $\eta_N := \sum_s x_s B^{(1)}_{\ell - nh,s}$ has the required properties.  This proves  the claim, and it gives  a proof of part (5) of the lemma. \B(5),\ref{L1}

\subsection{}{\bf Lemma}\label{L2} (the induction step). {\it Let $u \in \ZZ$ with $u \geq nT + 1$. Suppose $D_1 \in M^r$ such that $D_1 \equiv C_1 \pmod{(M_{1/T})^r}$, and such  that $\cF^h D_1 - p^n D_1 =: \xi \in (M_{u/T})^r$. Then there exists $\eta \in (M_{(u/T) - n})^r$ such that for $E_1 := D_1 - \eta$ we have  $\cF^h E_1 - p^n E_1  \in (M_{(u+1)/T})^r$ and  $E_1 \equiv C_1 \pmod{(M_{1/T})^r}$.}
\\
{\bf Proof.} We write $\xi = \xi_N + \xi_P$ according to $M = N \times P$. We conclude that $\xi_N \in (N \cap M_{u/T})^r$ and $\xi_P \in (P \cap (M_{u/T})^r$. Using \ref{L1}, (5), we construct $\eta_N \in (N \cap M_{1/T})^r$ \ \ such that $(\cF^h  - p^n) \eta_N \equiv \xi_N \pmod{(M_{(u+1)/T})^r}$.  As $M_{u/T} \subset M_n$ we can choose 
 $\eta_P := - p^{-n}\xi_P$; we have $\eta_P \in M_{(u/T) - n}^r \subset (M_{1/T})^r$. 
  With  $\eta:= \eta_N + \eta_P$ we see that 
$$(\cF^h - p^n) \eta \equiv \xi \pmod{(M_{(u+1)/T})^r} \quad\mbox{and}\quad \eta \in (M_{1/T})^r.$$ 
Hence $(\cF^h - p^n)(D_1 - \eta) \in (M_{(u+1)/T})^r$ and  we see that  $E_1 := D_1 - \eta$ has the required properties.  This proves the lemma.  

\B\ref{L2}

\subsection{}{\bf Proof of} \ref{filt}. {\bf (1)} {\it There exists $E_1 \in M^r$ such that $(\cF_n - p^n)E_1 = 0$ and $E_1 \equiv C_1 \pmod{(M_{1/T})^r}$}.\\
{\bf Proof.}  For $u \in \ZZ_{\geq nT+1}$ we write $D_1(u) \in M^r$ for a vector such that
$$D_1(u) \equiv C_1 \pmod{(M_{1/T}} \quad\mbox{ and } \quad \cF^h D_1(u) - p^n D_1(u)  \in (M_{u/T})^r.$$
By \ref{L1},  (4),  the vector $C_1 =: D_1(nT+1)$ satifies this condition for $u = nT+1$. Here we start induction. By repeated application of \ref{L2} we conclude there exists a sequence 
$$\{D_1(u) \mid u \in \ZZ_{\geq nT+1}\}
\quad\mbox{such that}\quad D_1(u) - D_1(u+1) \in  (M_{(u/T) - n})^r$$
satisfying the conditions above. As $\cap_{\rho \to \infty} \ \ M_{\rho}= \{0\}$ this sequence converges. Writing $E_1 := D_1(\infty)$ we achieve the conclusion. \B(1)

\vn
{\bf (2)} {\it For every $j \geq 0$ we have}
$$p^{-[\frac{jn}{h}]}\cF^jE_1 \in M \ \ \forall j \geq 0;
\quad \textrm{define}\quad  
N' := \prod_{1 \leq j < h} \ \ W{\cdot}p^{-[\frac{jn}{h}]}\cF^jE_1 \subset M.$$
{\it This is a Dieudonn\'e submodule. Moreover there is an isomorphism 
$$\DD((H_{m,n})^r) \cong N',$$
 $N' \prod P \to N' + P$ is an isomorphism of $W$-modules and $N' + P = M$.  This constructs $X_1 \subset X$, with  
$$\DD(X_1 \subset X) = (N' \subset M)\quad\textrm{such that}\quad (X/X_1)[p] \cong \prod_{j>1}(M_{m_j,n_j})^{r_j}.$$  
}\\
{\bf Proof.}
By \ref{L1},  (2), we see that $\cF^jE_1 \in M_{[jn/h]}$, hence the first statement follows.
 
As $\cF^hE_1 = p^nE_1$ it follows that $N' \subset M$ is a Dieudonn\'e submodule; using \ref{struct} this shows $\DD((H_{m,n})^r) \cong N'$.\\
{\bf Claim.} {\it The images  $N' \twoheadrightarrow N' \otimes k = N'/pN' \subset M/pM$ and  $P \twoheadrightarrow P/pP \subset M/pM$ inside $M/pM$ have zero intersection and  $N' \otimes k + P \otimes k = M/pM$.} Here we write $- \otimes k = - \otimes_W (W/pW)$.\\
For $y \in \ZZ_{\geq 0}$ we write $g(y) := yn - h{\cdot}[\frac{yn}{h}]$; note that, in the notation in \ref{notat3}, we have $$p^{-[\frac{jn}{h}]}\cF^jC'_1 = A^{(1)}_{g(j)}.$$ 
Suppose
$$\tau := \sum_{0 \leq j < h} \beta_{j,s}  p^{-[\frac{jn}{h}]}\cF^j{\cdot}(E_{1,s} \bmod pM) \in \left(N' \otimes k \cap  P \otimes k\right) \subset M/pM, \ \  \beta_j \in k$$
such that  $\tau \not=0$. Let $x, s$ be  a pair of indices such that $\beta:= \beta_{x,s} \not=0$ and for every $y$ with $g(y) < g(x)$  we have $\beta_{y,s} = 0$. Project inside $M/pM$ on the factor $N_s$.  Then  $$\tau_s \equiv \beta{\cdot}B^{(1)}_{g(x),s} \pmod{M_{\frac{g(x)}{h} + \frac{1}{T}}+P},$$
which is a contradiction with the fact that $N \cap P = 0$ and with the fact that the residue class of 
$$B^{(1)}_{g(x),s} \quad\mbox{generates}\quad   \left((M_{\frac{g(x)}{h} }+P) /( M_{\frac{g(x)}{h} + \frac{1}{T}}+P)\right)_s = N_{\frac{g(x)}{h},s} / N_{\frac{g(x)}{h} + \frac{1}{h},s}.$$
We see that $\tau \not= 0$ leads to a contradiction.
This shows that  $N' \otimes k \cap P \otimes k = 0$ and $N' \otimes k + P \otimes k = M/pM$. Hence the claim is proved.

As $(N' \cap P) \otimes k \subset N' \otimes k \cap P \otimes k = 0$ this shows $(N' \cap P) \otimes k =0$. By Nakayama's lemma this implies $N' \cap P = 0$. The proof of the remaining statements follows.  This finishes the proof of {\bf (2)}, and it ends the proof of  the proposition.  \B\ref{filt}

\section{Split extensions and proof of the theorem}\label{EP}
\n
In this section we prove a proposition on split extensions. We will see that  Theorem \ref{th}  follows.

\subsection{}\label{ext}{\bf Proposition.}   {\it Let $(m, n)$ and $(d, e)$ be ordered pairs of pairwise coprime positive integers. Suppose that $n/(m+n) < e/(d+e)$.  Let 
$$ 0 \to Z:=H_{m,n}   \longrightarrow T
 \longrightarrow Y:=H_{d,e} \to 0$$
be an exact sequence of $p$-divisible groups such that
the induced sequence of the $p$-kernels spits:
$$0 \rightarrow Z[p] \stackrel{\leftarrow}\longrightarrow T[p] \stackrel{\leftarrow}{\longrightarrow} Y[p]  \rightarrow  0.$$
Then the sequence of $p$-divisible groups splits: $T \cong Z \oplus Y$.}

\subsection{}{\bf Remark.}  It is easy to give examples of a non-split extension $T/Z \cong Y$ of $p$-divisible groups, with  $Z$ non-minimal or $Y$ non-minimal, such that $T[p]/Z[p] \cong Y[p]$  splits.

\subsection{Proof of} \ref{th}. The theorem follows from \ref{filt2} and \ref{ext}. \B\ref{th}

\subsection{}\label{<} {\it In order to show} \ref{ext} {\it it suffices to prove}  \ref{ext} {\it under the extra condition  that $\frac{1}{2} \leq e/(d+e)$.}\\
In fact, if $n/(m+n) < e/(d+e) < \frac{1}{2}$,  we consider the exact sequence
$$ 0 \to H_{d,e}^t = H_{e,d}    \longrightarrow T^t
 \longrightarrow H_{m,n}^t = H_{n,m} \to 0$$
with  $\frac{1}{2} < d/(e+d) < m/(n+m).$  \B\ref{<}

\vn
{}From now on we assume that $\frac{1}{2} \leq e/(d+e)$.

\subsection{}\label{notat4} We fix notation which will be used in the proof of \ref{ext}. We write the Dieudonn\'e modules as: $\DD(Z)=N$, $\DD(T)=M$ and $\DD(Y)=Q$; we obtain an exact sequence of Dieudonn\'e modules $M/N = Q$, which is a split  exact sequence of $W$-modules, where $W = W_{\infty}(k)$. We write $m+n = h$ and $d+e = g$. We know that $Q$ is generated by elements $A_i$, with $i \in \ZZ_{\geq 0}$ such that $\pi(A_i) = A_{i+1}$, where $\pi \in \End(Q)$ is the uniformizer, and $\cV{\cdot}A_i = A_{i+d}$, \ \ $\cF{\cdot}A_i = A_{i+e}$; we know that $\{A_i \mid 0 \leq i < g = d+e\}$ is a $W$-basis for $Q$. Because $\frac{1}{2} \leq e/(d+e)$, hence $e \geq d$ we can 
choose generators for the Dieudonn\'e module $Q$ in the following way. We  choose integers $\delta_i$ by:
$$i{\cdot}g \leq \delta_i{\cdot}d < (i+1){\cdot}d + i{\cdot}e = ig + d$$
and integers $\gamma_i$ such that: 
$$\delta_1 = \gamma_1 + 1, \cdots ,\delta_i = \gamma_1 + 1 + \gamma_2 + 1 + \cdots +\gamma_i + 1;$$
note that $\delta_d = g = d+e$.  
We choose $C= A_0 = C_1$ and $\{C_1, \cdots ,C_d\} = \{A_0, \cdots ,A_{d-1}\}$ such that:
$$\cV^{\gamma_i} C_i = \cF C_{i+1}, \ \ 1 \leq i < d,  \ \ \cV^{\gamma_d}C_d = \cF C_1,  \quad\mbox{hence}\quad \cV^{\delta_i} C = p^i{\cdot}C_{i+1}, \ \ \ 1 \leq i < d;$$
note that $\cV^g C = p^d{\cdot} C$. With these choices we see that 
$$\{p^{-[\frac{jd}{g}]}\cV^jC\mid\ 0 \leq j < g\} = \{\cV^j C_i \mid 1 \leq i \leq d, \ \ 0 \leq j \leq \gamma_i\} = \{A_{\ell} \mid 0 \leq \ell < g\}.$$
Choose an element $B = B_1 \in M$ such that
$$M \longrightarrow Q \quad\mbox{gives}\quad B_1 = B \mapsto \left(B \ {\rm mod}\  N\right) = C = C_1.$$
Let $\pi'$ be the uniformizer of $\End(N)$. Consider the filtration $N = N^{(0)} \supset \cdots \supset   N^{(i)}   \supset N^{(i+1)} \supset \cdots$ defined by $(\pi')^i(N^{(0)}) =  N^{(i)}$. Note that $\cF  N^{(i)} =   N^{(i+n)}$, and $\cV  N^{(i)} =  N^{(i+m)} $, and $p^iN = N^{(i{\cdot}h)}$ for $i \geq 0$.

\subsection{}{\bf Proof of \ref{ext}.}\\
{\bf (1)} Construction of $\{B_1, \cdots , B_d\}$.  {\it For every choice of $B = B_1 \in M$ with $\left(B \ {\rm mod}\  N\right) = C$, and every  $1 \leq i < d$ we claim that $\cV^{\delta_i}B$ is divisible by $p^{i}$. Defining $B_{i+1} := p^{-i}\cV^{\delta_i}B$, we see that $B_i \ {\rm mod} \  N = C_i$ for $1 \leq i \leq d$. Moreover we claim: } 
$$\cV^g B - p^d{\cdot}B \in N^{(dh + 1)}.$$
Choose $B''_i \in M$ with  $B''_i  \ {\rm mod} \  N = C_i$. Then $\cV^{\gamma_i}B''_i - \cF B''_{i+1} =: p{\cdot}\xi_i \in pN$; hence $\cV^{\gamma_i +1}B''_i - p{\cdot}B''_{i+1} = p\cV \xi_i \in p\cV N$. For $1 < i \leq d$ we obtain that 
$$\cV^{\delta_i} B - p^{i}{\cdot}B = \sum_{1 \leq j <i} \cV^{\delta_i - \delta_j}p^j\cV \xi_j, \ \ \ \xi_j \in N.$$
{}From $n/(m+n) < e/(d+e)$ we conclude $g/d > h/m$; using
$\delta_i{\cdot}d \geq ig$ and $\delta_jd < (j+1)d + je$ we see:
$$i>j \quad\mbox{implies}\quad \delta_i - \delta_j + 1 > (i-j)(g/d) > (i-j)(h/m);$$
hence 
$$(\delta_i - \delta_j)m + j(m+n) + m > ih;$$
This shows 
$$\cV^{\delta_i - \delta_j}p^j\cV \xi_j \in p^iN^{(1)}.$$  
As $\delta_d = g$ we see that $\cV^g B - p^d{\cdot}B \in p^dN^{(1)} = N^{(dh + 1)}.$ 
\B(1)

\vn
{\bf (2)} The induction step. {\it Suppose that for a choice $B \in M$ with $\left(B \ {\rm mod}\  N\right) = C$, there exists an integer $s \geq dh + 1$ such that $\cV^g B - p^d{\cdot}B \in N^{(s)}$; then there exists a choice $B' \in M$  such that $B' - B \in N^{(s-dh)}$  and}
$$\cV^g B' - p^d{\cdot}B' \in N^{(s+1)}.$$
In fact, write  $p^d{\cdot}B - \cV^g B   = p^d{\cdot}\xi$. Then $\xi \in N^{(s-dh)}$. Choose $B' := B - \xi$. Then:
$$\cV^g B' - p^d{\cdot}B' = \cV^g B - p^d{\cdot}B - \cV^g \xi + p^d \xi = -\cV^g \xi \in N^{(gm-dh+s)};$$
note that $gm - dh >0$.
\B(2)

\vn
{\bf (3)} {\it  For any integer $r \geq d+1$, and $w \geq rh$
there exists $B=B_1$ as in} \ref{notat4} {\it  such that $\cV^g B - p^d B \in N^{(w)} = p^r{\cdot}N^{(w-rh)}$. This gives a homomorphism $\va_{r-d}$}
$$M/p^{r-d}M \longleftarrow Q/p^{r-d}Q
\quad\mbox{extending}\quad   M/pM \longleftarrow Q/pQ.$$  
The induction step (1) proves the first statement, induction starting at $w = dh+1$. Having chosen $B_1$, using (2) we construct $B_{i+1} := p^{-i}\cV^{\delta_i} B_1$ for $1 \leq i < d$. In that case on the one hand $\cV^{\gamma_d}B_d - \cF B_1 = p{\cdot}\xi_d$, on the other hand  $\cV^g B - p^d B \in N^{(w)} \subset p^rN$. Hence $p^d\cV \xi_d \in p^rN$; hence $p\xi_d \in p^{r-d}N$. This shows that the residue classes of $B_1, \cdots , B_d$ in  $M/p^{r-d}M$ generate a Dieudonn\'e module isomorphic to $Q/p^{r-d}Q$ which moreover by \ref{notat4} extends the given isomorphism induced by the splitting. \B(3)

\vn
By \cite{FO-Foliations}, 1.6 we see that for some large $r$ the existence of $M/p^{r-d}M \longleftarrow Q/p^{r-d}Q$ as in (3) shows that its restriction $M/pM \longleftarrow Q/pQ$ lifts to a homomorphism $\va$ of Dieudonn\'e modules $M \leftarrow Q$; in that case $\va_1$ is injective. Hence $\va$ splits the extension $M/N \cong Q$. 
Taking into account \ref{<} this   proves the proposition. \B\ref{ext}

\vn
{\bf Remark.} Instead of the last step of the proof above, we could construct an infinite sequence $\{B(u)\mid u \in \ZZ_{(d+1)h}\}$ such that  $\cV^gB(u) - p^dB \in N^{(u)}$  and $B(u+1)-B(u) \in N^{(u-dh)}$ for all $u \geq (d+1)h$. This sequence converges and its limit $B(\infty)$ can be used to define the required section.

\section{Some comments} 
\subsection{}\label{Rk41}{\bf Remark.} For any $G$,  a \BT1 over $k$, which is  {\it not minimal} there exist infinitely many mutually non-isomorphic $p$-divisible groups $X$ over $k$ such that $X[p] \cong G$. A central leaf and an EO-stratum are equal if and only if we are in the minimal situation. Details will appear in a later publication, see \cite{FO-Simple}.

\subsection{}{\bf Remark.} Suppose that $G$ is a minimal \BT1; we can recover the Newton polygon $\beta$ with the property $H(\beta)[p] \cong G$ from $G$. This follows from the theorem, but there are also other ways to prove this fact.

\subsection{} For \BT1 group schemes we can define a Newton polygon; let $G$ be a \BT1 group scheme over $k$, and let $G = \oplus G_i$ be a decomposition into indecomposable ones, see \cite{Kraft}. Let $G_i$ be of rank $p^{h_i}$, and let $n_i$ be the dimension of the tangent space of $G_i^D$; define $\cN'(G_i)$ be the isoclinic polygon consisting of $h_i$ slopes equal to $n_i/h_i$;  arranging the slopes in non-decreasing order, we have defined $\cN'(G)$. For a $p$-divisible group $X$ we compare $\cN(X)$ and $\cN'(X[p])$; these polygons have the same endpoints; some rules seem to apply; if $X$ is minimal, equivalently $X[p]$ is minimal, then $\cN(X) = \cN'(X[p])$. Besides this I do not see rules describing the relation between $\cN(X)$ and $\cN'(X[p])$. For Newton polygons $\beta$ and $\gamma$ with the same end points we write $\beta \prec \gamma$ if every point of $\beta$ is on or below $\gamma$. Note:
\begin{itemize}
\item There exists a $p$-divisible group $X$ such that $\cN(X) \succneqq \cN'(X[p])$; indeed, choose $X$ isosimple, hence  $\cN(X)$ isoclinic, such that $X[p]$ is decomposable.
\item There exists a $p$-divisible group $X$ such that $\cN(X) \precneqq \cN'(X[p])$; indeed, choose $X$ such that  $\cN(X)$ is not isoclinic, hence $X$ not isosimple,  all slopes strictly between $0$ and $1$ and $a(X) = 1$; then $X[p]$ is indecomposable, hence $\cN'(X[p])$ is isoclinic.
\end{itemize}
It could be useful to have better insight in the relation between various properties of $X$ and $X[p]$.


\blankline
\noindent
\begin{tabbing}
\hspace{50 mm}
        \= \= \kill 
        \=Frans Oort \> \\ 
        \>Mathematisch Instituut \> \\
        \>Budapestlaan 6 \> Postbus 80010\\   
       \>NL - 3584 CD TA Utrecht\> NL - 3508 TA Utrecht\\ 
       \>The Netherlands\> The Netherlands\\ 
        \>email: oort@math.uu.nl  \\
\end{tabbing}


\end{document}